\documentclass{amsart}

\usepackage{fancyhdr}
\usepackage{lastpage}
\usepackage{stmaryrd,yhmath}

\newcommand{\bdis}{\begin{displaymath}}
\newcommand{\edis}{\end{displaymath}}
\newcommand{\be}{\begin{equation}}
\newcommand{\ee}{\end{equation}}
\newcommand{\mbb}{\mathbb}
\newcommand{\mcal}{\mathcal}

\newcommand{\vp}{\varphi}

\newcommand{\mT}{\mathring{T}}

\newcommand{\zf}{\zeta\left(\frac{1}{2}+it\right)}

\newtheorem{theorem}{Theorem}

\theoremstyle{definition}

\newtheorem{cor}[]{Corollary}

\theoremstyle{remark}
\newtheorem{remark}[]{Remark}

\numberwithin{equation}{section}

\begin{document}

\title{Jacob's ladders, Bessel's functions and the asymptotic solutions of a new class of nonlinear integral equations}

\author{Jan Moser}

\address{Department of Mathematical Analysis and Numerical Mathematics, Comenius University, Mlynska Dolina M105, 842 48 Bratislava, SLOVAKIA}

\email{jan.mozer@fmph.uniba.sk}


\begin{abstract}
It is shown in this paper that there is a connection between the Riemann zeta-function $\zf$ and the Bessel's functions. In this direction, a new class 
of the nonlinear integral equations is introduced. 
\end{abstract}

\maketitle

\section{The first result} 

\subsection{}  

We obtain some new properties of the signal 
\bdis 
Z(t)=e^{i\vartheta(t)}\zf 
\edis  
that is generated by the Riemann zeta-function, where 
\bdis 
\vartheta(t)=-\frac t2\ln\pi+\text{Im}\ln\Gamma\left(\frac 14+i\frac t2\right)=\frac t2\ln\frac{t}{2\pi}-\frac t2-\frac{\pi}{8}+\mcal{O}\left(\frac 1t\right) , 
\edis 
namely, the properties connected with the interaction of the function $\zf$ with the Bessel's functions 
\bdis 
J_\nu(x)=\sum_{r=0}^\infty \frac{(-1)^r}{r!\Gamma(\nu+r+1)}\left(\frac{x}{2}\right)^{\nu+2r} 
\edis  
where $x>0,\ \nu>-1$ (this is the sufficient case for our pourpose). Let us remind that 
\bdis 
\tilde{Z}^2(t)=\frac{{\rm d}\vp_1(t)}{{\rm d}t},\ \vp_1(t)=\frac 12\vp(t) 
\edis  
where 
\be \label{1.1} 
\tilde{Z}^2(t)=\frac{Z^2(t)}{2\Phi^\prime_\vp[\vp(t)]}=\frac{Z^2(t)}{\left\{ 1+\mcal{O}\left(\frac{\ln\ln t}{\ln t}\right)\right\}\ln t} , 
\ee  
(see \cite{1}, (3.9); \cite{3}, (1.3); \cite{7}, (1.1), (3.1), (3.2)), and $\vp(t)$ is the Jacob's ladder, i.e. a solution to the nonlinear integral 
equation (see \cite{1}) 
\bdis 
\int_0^{\mu[x(T)]}Z^2(t)e^{-\frac{2}{x(T)}t}{\rm d}t=\int_0^TZ^2(t){\rm d}t . 
\edis  

\subsection{} 

The system of the Bessel's functions 
\bdis 
\{ J_\nu(\mu_n^{(\nu)}x)\}_{n=1}^\infty,\ x\in [0,1],\ J_\nu(\mu_n^{(\nu)})=0 
\edis 
is the orthogonal system on the segment $[0,1]$ with the weight $x$, i.e. the following formulae hold true 
\be \label{1.2} 
\begin{split} 
& \int_0^1J_\nu(\mu_m^{(\nu)}x)J_\nu(\mu_n^{(\nu)}x)x{\rm d}x=0 ,\ m\not= n , \\ 
& \int_0^1\left[ J_\nu(\mu_n^{(\nu)}x)\right]^2x{\rm d}x=\frac{1}{2}J_{\nu+1}(\mu_n^{(\nu}) . 
\end{split}
\ee 
It is shown in this paper that the $\tilde{Z}^2$-transformation of the Bessel's functions generates a new orthogonal system of functions that is 
connected with $\left|\zf\right|^2$. Namely, the following theorem holds true. 

\begin{theorem}  
Let $x=t-T,\ t\in [T,T+1]$ and 
\bdis 
\vp_1\{[\mT,\widering{T+1}]\}=[T,T+1],\ T\geq T_0[\vp_1] . 
\edis 
Then the system of functions 
\bdis 
J_\nu[\mu_n^{(\nu)}(\vp_1(t)-T)],\ t\in [\mT,\widering{T+1}],\ n=1,2,\dots 
\edis 
is orthogonal on $[\mT,\widering{T+1}]$ with the weight 
\bdis 
(\vp_1(t)-T)\tilde{Z}^2(t) , 
\edis  
i.e. the following system of the new-type integrals 
\be \label{1.3} 
\begin{split} 
& \int_{\mT}^{\widering{T+1}}J_\nu[\mu_m^{(\nu)}(\vp_1(t)-T)]J_\nu[\mu_n^{(\nu)}(\vp_1(t)-T)]\cdot \\ 
& \cdot (\vp_1(t)-T)\tilde{Z}^2(t){\rm d}t=0,\ m\not=n , \\ 
& \int_{\mT}^{\widering{T+1}}\{J_\nu[\mu_n^{(\nu)}(\vp_1(t)-T)]\}^2(\vp_1(t)-T)\tilde{Z}^2(t){\rm d}t=\frac 12[J_{\nu+1}(\mu_n^{(\nu)})]^2 
\end{split} 
\ee  
is obtained, where $\vp_1(t)-T \in  [0,1]$, and 
\be \label{1.4} 
\rho\{[0,1];[\mT,\widering{T+1}]\}\sim T, \ T\to \infty , 
\ee  
(where $\rho$ stands for the distance of corresponding segments). 
\end{theorem} 

\begin{remark} 
Theorem 1 gives the contact point between the functions $\zf,\ \vp_1(t)$  and the Bessel's functions $J_\nu(x)$. 
\end{remark} 

This paper is the continuation of the series \cite{1}-\cite{18}. 

\section{The second result: new class of nonlinear integral equations} 

\subsection{} 

Let us remind that 
\be \label{2.1} 
t-\vp_1(t)\sim (1-c)\pi(t) \ \Rightarrow \ \mT\sim T,\ T\to\infty 
\ee  
where $c$ is the Euler's constant and $\pi(t)$ is the prime-counting function. Then the second formula in (\ref{1.3}) via the mean-value theorem 
(comp. (\ref{1.1}) leads to 

\begin{cor} 
\be \label{2.2} 
\begin{split}  
& \int_{\vp_1^{-1}(T)}^{\vp_1^{-1}(T+1)}\{J_\nu[\mu_n^{(\nu)}(\vp_1(t)-T)]\}^2(\vp_1(t)-T)\left|\zf\right|^2{\rm d}t\sim \\ 
& \sim\frac 12[J_{\nu+1}(\mu_n^{(\nu)})]^2\ln T,\ T\to\infty,\ n=1,2,\dots \ . 
\end{split} 
\ee 
\end{cor} 

\begin{remark} 
Let the primary oscillations 
\bdis 
\left|\zf\right|,\ t\in[\vp_1^{-1}(T),\vp_1^{-1}(T+1)]
\edis  
interact with the complicated modulated oscillations 
\be \label{2.3} 
|J_\nu[\mu_n^{(\nu)}(\vp_1(t)-T)]|\sqrt{\vp_1(t)-T} . 
\ee  
Then the integral (\ref{2.2}) expresses the energy of the resulting oscillations. Let us note that the oscillations (\ref{2.3}) comers to the point 
$t$ with the big retardation (see (\ref{2.1})) 
\bdis 
t-\{\vp_1(t)-T\}=t-\vp_1(t)+T\sim (1-c)\pi(t)+T,\ T\to\infty . 
\edis 
\end{remark} 

\subsection{} 

In this direction the following theorem holds true. 

\begin{theorem} 
Every Jacob's ladder $\vp_1(t)=\frac 12\vp(t)$ where $\vp(t)$ is the (exact) solution to the nonlinear integral equation 
\bdis 
\int_0^{\mu[x(T)]}Z^2(t)e^{-\frac{2}{x(T)}t}{\rm d}t=\int_0^TZ^2(t){\rm d}t 
\edis  
is the asymptotic solution of the new-type nonlinear integral equation 
\be \label{2.4} 
\begin{split} 
& \int_{x^{-1}(T)}^{x^{-1}(T+1)}\{J_\nu[\mu_n^{(\nu)}(x(t)-T)]\}^2(x(t)-T)\left|\zf\right|^2{\rm d}t= \\ 
& =\frac{1}{2}[J_{\nu+1}(\mu_n^{(\nu)})]^2\ln T,\ T\to\infty,\ n=1,2,\dots \ . 
\end{split}
\ee 
At the same time the Jacob's ladder $\vp_1(t)$ is the asymptotic solution to the following nonlinear integral equations (comp. (\ref{2.2}) with 
\cite{18}, (1.5), (2.2), (3.2), (3.3), (3.5), (3.6)) 
\be \label{2.5} 
\begin{split} 
& \int_{x^{-1}(T)}^{x^{-1}(T+2)}[P_n^{\alpha,\beta}(x(t)-T-1)]^2(T+2-x(t))^\alpha(x(t)-T)^\beta\left|\zf\right|^2{\rm d}t= \\ 
& =\frac{2^{\alpha+\beta+1}}{2n+\alpha+\beta+1}\frac{\Gamma(n+\alpha+1)\Gamma(n+\beta+1)}{n!\Gamma(n+\alpha+\beta+1)}\ln T,\ n=1,2,\dots , 
\end{split}
\ee 

\be \label{2.6} 
\int_{x^{-1}(T)}^{x^{-1}(T+2)}[P_n(x(t)-T-1)]^2\left|\zf\right|^2{\rm d}t=\frac{2}{2n+1}\ln T,\ n=1,2,\dots , 
\ee 

\be \label{2.7} 
\int_{x^{-1}(T)}^{x^{-1}(T+2)}[T_n(x(t)-T-1)]^2\frac{\left|\zf\right|^2}{\sqrt{1-(x(t)-T-1)^2}}{\rm d}t=\frac{\pi}{2}\ln T,\ n=1,2,\dots , 
\ee 

\be \label{2.8} 
\int_{x^{-1}(T)}^{x^{-1}(T+2)}\frac{\left|\zf\right|^2}{\sqrt{1-(x(t)-T-1)^2}}{\rm d}t=\pi\ln T , 
\ee 

\be \label{2.9} 
\begin{split} 
& \int_{x^{-1}(T)}^{x^{-1}(T+2)}[U_n(x(t)-T-1)]^2\sqrt{1-(x(t)-T-1)^2}\left|\zf\right|^2{\rm d}t= \\ 
& = \frac{\pi}{2}\ln T,\ n=1,2,\dots , 
\end{split} 
\ee  

\be \label{2.10} 
\int_{x^{-1}(T)}^{x^{-1}(T+2)}\sqrt{1-(x(t)-T-1)^2}\left|\zf\right|^2{\rm d}t=\frac{\pi}{2}\ln T , 
\ee 

for every fixed $T\geq T_0[\vp_1]$ where $P^{\alpha,\beta}_n(t),\ P_n(t),\ T_n(t),\ U_n(t)$ denote the polynomials of Jacobi, Legendre and 
Chebyshev of the first and second kind, respectivelly. 
\end{theorem} 

\begin{remark} 
There are the fixed-point methods and other methods of the functional analysis used to study the nonlinear equations. What can be obtained by using these 
methods in the case of the nonlinear integral equations (\ref{2.4})-(\ref{2.10})? 
\end{remark} 

\section{Proof of Theorem 1} 

\subsection{} 

Let us remind that the following lemma holds true (see \cite{6}, (2.5); \cite{7}, (3.3)): for every integrable function (in the Lebesgue sense) 
$f(x),\ x\in [\vp_1(T),\vp_1(T+U)]$ we have 

\be \label{3.1} 
\int_T^{T+U}f[\vp_1(t)]\tilde{Z}^2(t){\rm d}t=\int_{\vp_1(T)}^{\vp_1(T+U)}f(x){\rm d}x,\ U\in \left.\left( 0,\frac{T}{\ln T}\right.\right] 
\ee  
where 
\bdis 
t-\vp_1(t)\sim (1-c)\pi(t) , 
\edis  
$c$ is the Euler's constant and $\pi(t)$ is the prime-counting function. In the case (comp. Theorem 1) $T=\vp_1(\mT), T+U=\vp_1(\widering{T+U})$ we 
obtain from (3.1) 
\be \label{3.2} 
\int_{\mT}^{\widering{T+U}}f[\vp_1(t)]\tilde{Z}^2(t){\rm d}t=\int_T^{T+U}f(x){\rm d}x . 
\ee  

\subsection{} 

Putting 
\bdis 
f(t)=J_\nu[\mu_m^{(\nu)}(t-T)]J_\nu[\mu_n^{(\nu)}(t-T)](t-T),\ U=1 
\edis  
we have by (\ref{3.2}) and (\ref{1.2}) the following $\tilde{Z}^2$-transformation 
\bdis 
\begin{split} 
& \int_{\mT}^{\widering{T+1}}J_\nu[\mu_m^{(\nu)}(\vp_1(t)-T)]J_\nu[\mu_n^{(\nu)}(\vp_1(t)-T)]\cdot \\ 
& \cdot (\vp_1(t)-T)\tilde{Z}^2(t){\rm d}t= \\ 
& =\int_{T}^{T+1}J_\nu[\mu_m^{(\nu)}(t-T)]J_\nu[\mu_n^{(\nu)}(t-T)](t-T)\tilde{Z}^2{t}{\rm d}t= \\ 
& = \int_0^1J_\nu(\mu_m^{(\nu)}x)J_\nu(\mu_n^{(\nu)}x)x{\rm d}x=0, \ m\not=n , 
\end{split}
\edis 
where $t=x+T$, i.e. the first formula in (\ref{1.3}) holds true. Similarly, we obtain the second formula in (\ref{1.3}). 

\subsection{} 

Next, for $\xi\in (\mT,\widering{T+1})$ we have (see (\ref{2.1}) and \cite{18}, (4.4)) 

\be \label{3.3} 
\ln\xi=\ln\mT+\mcal{O}\left(\frac{1}{\ln\mT}\right)=\ln T+\mcal{O}\left(\frac{1}{\ln T}\right) . 
\ee  

The property (\ref{3.3}) was used in (\ref{2.1}), \dots .

\thanks{I would like to thank Michal Demetrian for helping me with the electronic version of this work.}


\begin{thebibliography}{29}
\bibitem{16}
J. Moser, `Jacob's ladders and the nonlocal interaction of the function $Z^2(t)$ with the function $\tilde{Z}^2(t)$ on the distance $\sim (1-c)\pi(t)$
for the collections of disconnected sets`, (2010), arXiv: 1007.5147.
\end{thebibliography}
\end{document}